\documentclass[12pt]{article}

\usepackage{amsmath}
\usepackage{amsfonts}
\usepackage{amssymb}
\usepackage{amsthm}
\theoremstyle{plain}
\newtheorem{thm}{Theorem}[section]
\newtheorem{prop}[thm]{Proposition}
\newtheorem{defn}{Definition}[section]

\newcommand{\at}[2][]{#1|_{#2}}

\numberwithin{equation}{section}

\begin{document}

\title{Arithmetic Properties of Sparse Subsets of $\mathbb{Z}^n$}
\author{Paul Potgieter}
\date{}
\maketitle \vspace{-1.2cm}
\begin{center}
 \emph{Department of Decision Sciences, University
of South Africa\\ P.O. Box 392, Pretoria 0003, South Africa}\\
\texttt{potgip@unisa.ac.za}
\end{center}
\setcounter{section}{0}

\begin{abstract}

Arithmetic progressions of length three may be found in compact subsets of the reals that satisfy certain Fourier- as well as Hausdorff-dimensional requirements. Similar results hold in the integers under analogous conditions, with  Fourier dimension being replaced by the decay of a discrete Fourier transform. In this paper we make this correspondence more precise, using a well-known construction by Salem. Specifically, we show that a subset of the integers can be mapped to a compact subset of the continuum in a way which preserves certain dimensional properties as well as arithmetic progressions of arbitrary length. The higher-dimensional version of this construction is then used to show that certain parallelogram configurations must exist in sparse subsets of $\mathbb{Z}^n$ satisfying appropriate density and Fourier-decay conditions. 

\medskip



\end{abstract}

\section{Introduction}

Subsets of $\mathbb{R}^n$ satisfying certain Hausdorff- and Fourier-dimensional properties have been shown to have interesting arithmetic properties \cite{chan2016finite, LP, shmerkin2015}. To obtain sets with the relevant properties, it has been standard practice to use integer sequences to construct Cantor-type subsets of $\mathbb{R}^n$. However, little has been said about the interpretation of such results in subsets $\mathbb{Z}^n$. It seems likely that such interpretations should exist, because similar notions of dimension and of Fourier decay can be formulated in both the continuum and in $\mathbb{Z}^n$. Yet the translation of results from the continuum to the discrete case is not necessarily trivial and there is at present no canonical method for doing so. We have, for example, the following result in the continuum:

 \begin{thm}\cite{LP}\label{thm2.1}
Assume that $E\subset[0,1]$ is a closed set which supports a probability measure $\mu$
with the following properties:
\begin{enumerate}
\item[(i)]{$\mu([x,x+\epsilon])\leq C_1\epsilon^\alpha$ for all $0<\epsilon\leq 1$,}
\item[(ii)]{$|\hat{\mu}(k)| \leq C_2 (1 - \alpha)^{-B} |k|^{-\frac{\beta}{2}}$ for all $k \ne 0$,}
\end{enumerate}
where $0<\alpha<1$ and $2/3<\beta\leq 1$.  If $\alpha>1-\epsilon_0$, where
$\epsilon_0>0$ is a sufficiently small constant depending only
on $C_1,C_2,B,\beta$, then $E$ contains a non-trivial 3-term arithmetic progression.
\end{thm}
The first condition of the theorem ensures that the Hausdorff dimension of $E$ is at least $\alpha$. Theorem \ref{thm2.2} shows that very similar conditions apply in $\mathbb{N}$ to obtain arithmetic progressions of length $3$. We first define an integer analogue to Hausdorff dimension~\cite{potgieter2}. Throughout the paper, $[0,N)$ denotes the interval $\{0,1,2,\ldots ,N-1 \}$ in $\mathbb{N}$. 

\begin{defn}
We say that a set $A\subseteq \mathbb{N}$ has \emph{upper fractional
density} $\alpha$ if
\[ \limsup_{N\to \infty} \frac{|A\cap [0,N)|}{N^{\beta}}\]
is $\infty$ for any $\beta < \alpha$ and $0$ for any $\beta >
\alpha$. We indicate this by the notation $d^{*}_f (A) = \alpha$.
\end{defn}
The lower fractional density can be similarly defined, and if the upper and lower fractional densities coincide we can of course simply refer to the \emph{fractional density} of $A$, which will be denoted by $d_f (A)$. The definition can trivially be extended to all of $\mathbb{Z}$. We have the following analogue of Theorem \ref{thm2.1} in the integers, where now we use the discrete Fourier transform instead of the Fourier-Stieltjes transform of a measure. 

\begin{thm}\cite{potgieter2}\label{thm2.2}
Let $A\subseteq \mathbb{N}$. Suppose $A$ satisfies the following
conditions:
\begin{itemize}
\item[(i)]{$A$ has upper fractional density $\alpha$, where $\alpha >1/2$.}
\item[(ii)]{For some $C>0$, the Fourier coefficients of the indicator functions $\chi_{A_N}$ of $A_N = A\cap[0,N)$ satisfy
\[|\widehat{\chi_{A_N}} (k)| \leq C(kN)^{-\beta/2}\] for any $0<k<N$, for large $N$, for some $2/3< \beta \leq 1$ satisfying $\beta>2-2\alpha$.}
\end{itemize}
Then $A$ contains a non-trivial arithmetic progression of length 3.
\end{thm}

Examples of sets satisfying the conditions of Theorems \ref{thm2.1} and \ref{thm2.2} were constructed by \L{}aba and Pramanik \cite{LP} and Potgieter \cite{potgieter2}, respectively. One could ask whether the results could be extended to longer arithmetic progressions by, for instance, requiring a larger Hausdorff or Fourier dimension. This seems to be unlikely, since these dimensional conditions do not seem to provide enough structure to the set for longer progressions to be inevitable. Tao \cite{Tao} discusses in full the failure of Fourier-dimensional conditions to guarantee progressions larger than $3$. As was shown by Shmerkin \cite{shmerkin2017salem}, given equal Fourier and Hausdorff dimensions, the dimension can be increased arbitrarily whilst avoiding the existence of \emph{any} non-trivial arithmetic progressions. 

Theorem \ref{thm2.2} was proved independently of Theorem \ref{thm2.1}, but using related techniques. The question is whether there is some method which can extract results such as Theorem \ref{thm2.2} from results like Theorem \ref{thm2.1}. This paper partially answers this question in two cases, namely deducing the existence of arithmetic progressions in subsets of $\mathbb{Z}$ or parallelogram configurations in subsets of $\mathbb{Z}^n$ from their existence in corresponding subsets of $\mathbb{R}$ or $\mathbb{R}^n$, respectively. The techniques utilized are mostly elementary, and seem to indicate that arithmetic properties analytically obtained in $\mathbb{R}^n$ may be interpreted in $\mathbb{Z}^n$, as long as one is careful with the construction used.

The basis of these results is the examination of the construction of Salem-type sets in the continuum using integer sequences. These methods are not dissimilar to the construction of the standard triadic Cantor set, although the specific construction used in this paper is a version of one used by Salem himself~\cite{Salem}. They form the core of much of the study of sparse sets in $\mathbb{R}^n$, including important restriction phenomena~\cite{bluhm1998theorem,hambrook2013sharpness,Mockenhaupt}. The construction is presented here not only for completeness, but because it is an integral part of the results obtained. The construction for higher dimensions is also rarely presented in full in a manner which highlights the dependencies of chosen constants, and therefore including it promotes self-containment.

The notion of density that will be used throughout is Hausdorff dimension, and its analogue in the integers, fractional density. Uniformity will be defined in terms of the decay of a Fourier-Stieltjes transform of a measure on a compact subset of $\mathbb{R}^n$, and in $\mathbb{Z}^n$ by a decay condition on the Fourier transform of the indicator function of a set.  We define the Fourier-Stieltjes transform of a finite Borel measure $\mu$ supported on a compact set $A\subset \mathbb{R}^n$ as
\begin{equation}\label{eq1.1}
\hat{\mu}(\xi) = \int_A e^{-2\pi i\xi \cdot x}d\mu (x).
\end{equation}

The \emph{Fourier dimension} of a set is a measure of exactly how rapid the decay of the transform of that measure is. 

\begin{defn}\label{def1.2} A compact set $A \subseteq [0,1]^n$ is said to have Fourier dimension $\beta$ if \begin{equation*}
 \beta = \sup \{\alpha \in [0,n]: \exists M_{1}^{+} (A)\left( |\hat{\mu}(\xi) |^2 =  o(|\xi |^{-\alpha })\right) \}.
\end{equation*}
\end{defn}

Here $M_{1}^{+} (A)$ denotes the set of probability measures on the set $A$. Condition (ii) of Theorem \ref{thm2.1}, for example, ensures that the Fourier dimension of the set under consideration is at least $\beta$. Salem \cite{Salem1} first showed how to construct sets with a specific Fourier dimension, and in a way that it coincides with the set's Hausdorff dimension. Hence, sets of equal Fourier and Hausdorff dimension are called \emph{Salem sets}. It is possible for sets to have a positive Fourier dimension which is strictly smaller than their Hausdorff dimension (it can never be larger). Such sets are referred to as \emph{Salem-type} sets. 

Following Theorem \ref{thm2.2}, we can define the Fourier dimension of a subset of the positive integers.

\begin{defn}
The \emph{Fourier dimension} of set $A\subseteq \mathbb{N}$ is the supremum of all $\beta \in [0,1]$ for which there exist a constant $C>0$ such that, for $A_N = A\cap [0,N)$, 
\begin{equation*}
\frac{1}{N}\sum_{n=0}^{N-1} \chi_{A_N} (n) e^{-2\pi i \frac{mn}{N}} \leq C(mN)^{-\beta/2},
\end{equation*} 
for any $0<m<N$, $m\in \mathbb{N}$, for all large $N$.
\end{defn}
A Salem-type subset of the positive integers will then be one which has Fourier dimension strictly less than its fractional density.

In Section 2 we present the construction of a subset of $[0,1]$ from a given subset of the positive integers. Briefly, given some $A\subseteq \mathbb{N}$ and some strictly increasing sequence of positive integers $(N_i)_{i\geq 1}$, we form the subsets $A_i = A \cap [0, N_i )$. The interval $[0,1]$ is divided into $N_1$ equal intervals, and a subinterval $[i/k, i/k +\xi_1)$ is chosen whenever $i\in A_1$, where $0<\xi_1 <1/N_1$. The union of all subintervals chosen can be denoted by $B_1$. Once $B_k$ is chosen, we find $B_{k+1}$ by applying the construction applied to $[0,1]$ to each interval comprising $B_{k}$. That is, each interval in $B_k$ is divided into $N_{k+1}$ equal intervals, and from the $i$-th interval the initial subinterval of length $\xi_{k+1}$, $0<\xi_{k+1}<( N_1 \cdots N_{k+1})^{-1}$ is chosen whenever $i\in A_{k+1}$. The intersection of all the $B_k$, $k=1,2,\ldots$, is denoted by $B$, or $B(A)$ if we want to highlight the dependence on $A$. We then have the following proposition:    

\begin{prop}\label{prop2.3}
Suppose that $A\subseteq \mathbb{N}$ satisfies the conditions of Theorem \ref{thm2.2}. Then the set $B$ obtained from $A$ via the construction in Section \ref{sec2.2} has Fourier dimension no less than $\beta$.
\end{prop}

As is noted in Section 2, it is easily shown that the Hausdorff dimension of the set $B$ is equal to the upper fractional density of the set $A$. 

The fact that such dimensional concepts are preserved is an indication that arithmetic properties might also be preserved, as evinced by the correspondence between Theorems \ref{thm2.1} and \ref{thm2.2}. That this is indeed the case is illustrated by the following, which holds irrespective of Hausdorff and Fourier dimension.

\begin{prop}\label{prop2.4}
For any $n\in \mathbb{N}$, a set $A\subseteq \mathbb{N}$ contains an arithmetic progression of length $n$ if and only if the set $B(A)$ defined by \eqref{eq2.4} contains an arithmetic progression of length $n$.
\end{prop}

The construction can be generalized to higher dimensions in a straightforward way, as presented in Section \ref{sec3.2}, in a manner which preserves relevant dimensions similarly to the one-dimensional case. The concept of a Salem-type set in $\mathbb{Z}^n$ is defined analogously to the one-dimensional case (Definition \ref{def3.2}).  The type of linear configuration we will examine is given by the following.

\begin{defn}\label{def4.3}
Suppose that $n\geq 2$, $m\geq n$ and $k\geq 3$. Let $\mathbb{B} = \{ B_i$: $i=1,\ldots ,k\} $ be a collection of $(m-n) \times n$ matrices with non-negative integer entries. We say that a set $A\subseteq \mathbb{Z}^n $ contains a $\mathbb{B}$-configuration if there exist $x\in \mathbb{Z}^n$ and $y\in \mathbb{Z}^{m-n}\setminus \{0\}$ such that $x+B_i y \in A$ for each $i$. 
\end{defn}

As is shown in Section \ref{sec4.3}, we can define parallelograms in $\mathbb{Z}^n$ as such $\mathbb{B}$-configurations. Setting aside the technical definition of a non-degenerate collection of matrices until Section \ref{sec4.1}, the main result of the paper is as follows:

\begin{thm}\label{thm4.2}
Suppose
\begin{equation*}
n \left\lceil \frac{k+1}{2} \right\rceil = m 
\end{equation*}
and 
\begin{equation*}
\frac{2(nk-m)}{k} < \beta <n.
\end{equation*}
Let $\mathbb{B} = \{B_1 ,\ldots ,B_k \}$ be a collection of $n\times (m-n)$ matrices with non-negative integer entries such that $A_j = (I_{n\times n} B_j )$ is non-degenerate in the sense of Definition \ref{def4.2}. Furthermore, letting $b^{(l)}_{i,j}$ denote the entries of the matrix $B_l$, $l=1,\ldots ,k$, we require that 
\begin{enumerate}
\item[(i)]{$B_1$ is the $n\times (m-n)$ zero matrix,}
\item[(ii)]{for each $l = 1, \ldots , m/n$ and each $i=1,2,\ldots ,n$, $\sum_{j=1}^{n} b^{(l)}_{i,j}\leq 1$,}
\item[(iii)]{for each $l = m/n +1 ,\ldots ,k$, $B_l$ is of the form $B_i + B_j$ for some $i,j \in \{2,3,\ldots ,m/n\}$.}
\end{enumerate}
Then for any constant $C>0$ there exists a positive $\varepsilon = \varepsilon(C,n,k,m,\mathbb{B})\ll 1$ with the following property. Suppose that the set $A\subseteq \mathbb{Z}^n$ is of Salem-type, as in Definition \ref{def3.2}, with a fractional density $\alpha$, where $n-\varepsilon < \alpha <n$ and a Fourier dimension of at least $\beta$. Then $A$ contains a $\mathbb{B}$-configuration.
\end{thm}
 
The proof of Theorem \ref{thm4.2} is presented in Section \ref{sec4.2}. In Section \ref{sec4.3}, Theorem \ref{thm4.2} is applied to show the existence of parallelograms in certain subsets of $\mathbb{Z}^n$. Section \ref{sec4.4} discusses some possible generalizations of the current work.

\section{Correspondence between $\mathbb{Z}$ and $[0,1]$}

\subsection{Constructing sets in $[0,1]$ from the integers}\label{sec2.2}

The construction used in this section to generate subsets of $[0,1]$ from subsets of $\mathbb{Z}$ is essentially due to Salem~\cite{Salem}.  
We now proceed with the construction, adjusted to our needs. Let $A$ be a subset of the non-negative integers, with a fractional density of $\alpha>0$. Let $(\alpha_i )_{i\geq 1}$ be any real sequence strictly increasing to $\alpha$, and let $(\gamma_i )_{i\geq 1}$ be any real, increasing sequence tending to infinity. By the definition of upper density, we can find a sequence $(N_i )_{i\geq 1}$ (to be used in \eqref{eq2.9.1}) such that 
\begin{equation}\label{eq2.3}
c_i := \frac{|A\cap [0, N_i )|}{N_{i}^{\alpha_i}}>\gamma_i.
\end{equation}
Thus, $c_i \to \infty$.  
We can also observe that, given the fractional density of $A$, certainly
\begin{equation*}
|A\cap [0,N_i)| = o(2^{\beta} N_{i}^{\beta}) 
\end{equation*}
for any $\beta > \alpha$, and therefore that 
\begin{equation}\label{c_decay}
c_i = o(2^{\beta} N_{i}^{\beta - \alpha_i})
\end{equation}
for any $\beta > \alpha$. This observation is useful in the determination of the Hausdorff dimension of the set to be constructed.

Divide $[0,1]$ into $N_1$ equal intervals and choose an interval $[i/N_1, i/N_1 +\zeta_1)$, where $0< \zeta_1 < N_{1}^{-1}$, if $i\in A_1 = A\cap [0, N_1 )$. The set of intervals chosen is denoted by $B_1$. That is,
\begin{equation*}
B_1 = \bigcup_{i\in A_1} \left[ \frac{i}{N_1},\frac{i}{N_1}+\zeta_1 \right).
\end{equation*}
Therefore, $B_1$ consists of $|A_1 |$ intervals of length $\zeta_1$. Let $A_k = A\cap [0, N_k )$ and inductively define $B_{k+1}$ by dividing the intervals belonging to $B_k$ into $N_{k+1}$ equal intervals and setting 
\begin{equation*} \left[ x + \frac{\zeta_1 \cdots \zeta_k i}{N_{k+1}}, x + \frac{\zeta_1 \cdots \zeta_k i}{N_{k+1}}+ \zeta_1 \cdots \zeta_{k+1} \right) \subset B_{k+1}
\end{equation*} if 
\begin{equation}\label{eq2.3.1} \left[ x, x +\zeta_1 \cdots \zeta_k \right) \subseteq B_k
\end{equation} 
and $i\in A_{k+1}$, where $0< \zeta_{k+1} < N_{k+1}^{-1}$ and $x$ is of the form
\begin{equation*}
x = \frac{j_1}{N_1}+ \frac{\zeta_1 j_2}{N_2}+\cdots + \frac{\zeta_1 \cdots \zeta_{k-1} j_k}{N_k} 
\end{equation*}
for some $j_1\in A_1, \ldots ,j_k \in A_k$. Thus, the set $B_{k+1}$ consists of $|A_1 | \cdots |A_{k+1} |$ intervals, each of length $\zeta_1 \cdots \zeta_{k+1}$. More concisely by letting $J_k$ be the set of all $x$ such that \eqref{eq2.3.1} holds,
\begin{equation*}
B_{k+1} = \bigcup_{x\in J_k}\bigcup_{i\in A_{k+1}}\left[ x +\frac{\zeta_1 \cdots \zeta_{k} i}{N_{k+1}}, x +\frac{\zeta_1 \cdots \zeta_{k} i}{N_{k+1}}+\zeta_1 \cdots \zeta_{k+1}\right).
\end{equation*}
We define the set $B$ as 
\begin{equation}\label{eq2.4}
B=\bigcap_{k=1}^{\infty} B_k.
\end{equation} 

The set $B$ obtained in this manner is denoted by $B(\zeta_1, \zeta_2, \ldots ,\zeta_k ,\ldots)$, although we usually do not find it necessary to explicitly indicate the parameters used in the construction. If it is necessary to reference the set $A$ used in the construction, we denote the resulting set by $B(A)$.

In Salem's original construction, it was required that, if we set $\xi_i =  N_{i}^{-1}$ for each $i$, the interval lengths $\zeta_i$ satisfy
\begin{equation*}
 \left( 1- \frac{1}{2i^2 }\right) \xi_i \leq \zeta_i < \xi_i, \quad i=1,2,\ldots.
\end{equation*}
This condition is specific to Salem's approach and not necessary in our construction, and would in fact be detrimental to our goal of showing that arithmetic progressions are preserved by the construction. Thus, we set our parameters 
\begin{equation}\label{eq2.5.5}
\zeta_i = 1/2N_{i}
\end{equation}
 at each stage. To verify that the set $B$ has Hausdorff dimension equal to $\alpha$, the upper fractional density of $A$, is easy to do, since the fractional density of the original subset of the integers dictates exactly how many intervals there are at each stage of the construction, and the diameter of each is known. What is more, no interval is contiguous to any other, which makes the sum involved in determining Hausdorff dimension trivial to compute, guaranteeing a dimension of no less than $\alpha$. The condition \eqref{c_decay} in turn ensures that the Hausdorff dimension is no larger than $\alpha$, hence the desired result. (See, for instance, Example 4.6 of Falconer~\cite{falconer2004fractal}.)

There is a way of representing each element of the set $B$ as an infinite sum, supposing that $B$ was constructed as above. Let $a^{(k)}(i)$ denote the $i$-th member of the set $A_k$. Each point $x$ belonging to the set $B$ can then be written as
\begin{equation}\label{eq2.6}
 x = \frac{a^{(1)}(\varepsilon_1 )}{N_1 } + \frac{\zeta_1 a^{(2)}(\varepsilon_2)}{ N_2 }+\frac{\zeta_1 \zeta_2  a^{(3)}({\varepsilon_3})}{ N_3}+\cdots ,
 \end{equation} where $\varepsilon_i$ ranges over all values $1,2, \ldots, |A_i|$. This representation is useful in calculating Fourier-Stieltjes transforms of measures on the set $B$.

\subsection{Fourier correspondence between subsets of $\mathbb{N}$ and $[0,1]$}\label{sec2.3}

Since, as shown in the previous section, we can easily move from $\mathbb{N}$ to $[0,1]$ in a way which preserves a notion of Hausdorff dimension, the question becomes whether we can do the same with Fourier dimension. The key to this is to use the decay of the Fourier coefficients of the indicator functions of the sets $A_N$, seen as contained in the group $\mathbb{Z}_N$, to obtain bounds on the Fourier-Stieltjes transforms of certain discrete probability measures on $[0,1]$. These Fourier-Stieltjes transforms can be used to find the Fourier-Stieljes transform of a certain probability measure $\mu$ on the set $B$ constructed above.  

Letting $\chi_{A_N}$ denote the indicator function of $A_N$ considered as a subset of the cyclic group $\mathbb{Z}_N$, the discrete Fourier transform is defined as follows:
\begin{equation*}
\widehat{\chi}_{A_N} (k) = \frac{1}{N}\sum_{n=0}^{N-1} \chi_{A_N}(n) e^{-2\pi i\frac{kn}{N}}.
\end{equation*}

We remark that the construction of the measures in this section may seem slightly more deterministic in character than those constructed by Salem. His probabilistic construction yields the required dimension only almost surely, whereas our adapted construction will certainly yield a measure with the required decay. It has to be admitted that this only hides the probabilistic aspect, since the most likely way to obtain sets with the given properties would be through probabilistic means. A probabilistic construction of such a set is given in Potgieter~\cite{potgieter2}.

We now present the proof of Proposition \ref{prop2.3}. Although we cannot exclude the possibility that Fourier dimension will be preserved for smaller values, our proof does require the Fourier decay of the indicator function of $A$ to be rapid, and that it satisfies the above relation to the fractional density. This is sufficient for our purposes, because the configurations we are interested in seem to only occur under such conditions, even though weaker conditions cannot be entirely ruled out in exceptional circumstances.  

\begin{proof}[Proof of Proposition \ref{prop2.3}.] We construct a measure on the set $B$ obtained in the previous section, assuming the set $A$ and the numbers $N_i$, $\alpha_i$, $c_i$, $i\in \mathbb{N}$, to be the same as before. At the $k$-th stage of the construction we introduce a continuous, non-decreasing function $F_k :[0,1] \to [0,1]$, where $F_k(0)=0$,  with the property that it increases linearly by $(c_1 \cdots c_k N_{1}^{\alpha_1}\cdots N_{k}^{\alpha_k})^{-1}$ over each interval of length  $\zeta_1 \cdots \zeta_k$ comprising $B_k$, with all the constants in question satisfying the requirements \eqref{eq2.3}, \eqref{c_decay} and \eqref{eq2.5.5} of Section \ref{sec2.2}. 
On intervals not part of the set $B_k$, the function remains constant, equal to its previously attained value on $B_k$. Hence, for each $k$ we have a piecewise linear function increasing from $0$ to $1$ over the unit interval, which is a distribution function for a probability measure $\mu_k$ supported on the set $B_k$. These functions converge pointwise to a function $F$ which will be continuous and non-decreasing, with $F(0)=0$ and $F(1)=1$, as can be easily seen by showing they form a Cauchy sequence with respect to the uniform norm. The measure corresponding to $F$ will be denoted by $\mu$.

Our interest now lies in calculating asymptotic properties of the Fourier-Stieltjes transform of $\mu$ as defined in \eqref{eq1.1}. To do so, we must find a way of representing upper bounds on $|\hat{\mu}|$ in terms of the upper bounds of the Fourier coefficients of the sets $A_{N_i}$.

We let $|A_k |$ denote the number of intervals chosen in the $k$-th step of the construction of the set $B=B(\zeta_1 ,\zeta_2 ,\ldots )$ from the given set $A$. The numbers $a^{(k)}(j)$, $k,j=1,2,3,\ldots$, are as determined in \eqref{eq2.6}.
We define 
\begin{equation*}
Q^{(k)}(u) = \frac{1}{|A_k |}\sum_{j=1}^{|A_k |} e^{-2\pi iu \frac{a^{(k)}(j)}{N_k}},
\end{equation*}
for all $k\geq 1$. These may be seen as the Fourier-Stieltjes transforms of discrete measures on the unit interval. As in Salem~\cite{Salem}, the Fourier-Stieltjes transform of the measure $\mu$ can then be written as
\begin{equation*}
\hat{\mu}(u) = \lim_{m\to \infty}Q^{(1)}(u)\prod_{k=1}^{m}Q^{(k+1)}(\zeta_1 \cdots \zeta_k u)
\end{equation*}
(see also, for instance, p.\ 195 of Zygmund~\cite{Zygmund}).

Note that $Q^{(k)}$ corresponds closely to the discrete Fourier transform (modulo $N_k$) on integer values:

\begin{equation*}
Q^{(k)}(m) = \frac{1}{|A_{k}|}\sum_{n=0}^{N_k -1} \chi_{A_{k}} (n)e^{-2\pi in\frac{m}{N_k}}=
\frac{N_k}{|A_k |}\widehat{\chi_{A_{k}}}(m).
\end{equation*}

We can therefore dictate the decay of $\hat{\mu}$ by requiring sufficient decay of $\widehat{\chi_{A_{k}}}(m)$, as in the statement of the theorem. Thus, for any $\beta$ satisfying the conditions of the theorem, there exists some constant $C$ such that for any $m\in [1,N_k)\subset \mathbb{N}$,
\begin{equation}\label{eq2.8}
|Q^{(k)}(m)| \leq C\frac{N_k}{|A_k |}(mN_k)^{-\beta /2} = c_{k}^{-1} C m^{-\beta /2} N_{k}^{1-\frac{\beta}{2} - \alpha_k }.
\end{equation}
We will need to consider arbitrarily large arguments in the approximation of $|\hat{\mu}(u)|$, in which case the $m$ in the above expression ought to be replaced by $m\bmod N_k$. However, once $m\geq N_k$, we will use some other $|Q^{(j)}(m)|$ with $j>k$ as the main component of our approximation, thus avoiding the issue. This also lets us avoid the problem of approximating $|Q^{(k)}(m)|$ when $m$ is a multiple of $N_k$, where it would be equal to $1$. 

Since the inequalities bounding $|Q^{(k)}|$ only hold for integer-valued arguments, we will need to extend them to arbitrary real arguments in order to approximate $|\hat{\mu}|$. To do so, note that an elementary calculation gives the upper bound
\begin{equation*}
\left| \frac{d}{du}|Q^{(k)}(u)|  \right|  \leq |Q^{(k)}(u)|
\end{equation*}
for any $u>0$. For any interval $[u,u+a]$, $a>0$, we then have that 
\begin{align*}
\max_{x\in [0,a]} |Q^{(k)}(u+x)| &\leq a\max_{x\in [0,a]} \left| \frac{d}{dw}|Q^{(k)}(w)| \at[\big]{w=u+x} \right| +|Q^{(k)}(u)| \\
&\leq  a \max_{x\in [0,a]} |Q^{(k)}(u+x)| +|Q^{(k)}(u)|
\end{align*}
and therefore 
\begin{equation*}
\max_{x\in [0,a]} |Q^{(k)}(u+x)| \leq \frac{1}{1-a}|Q^{(k)}(u)| 
\end{equation*}
when $1-a>0$. Set $a=1/2$. For any integer $m$ we therefore have that for any $u\in [m,m+1/2]$,
\begin{equation*}
|Q^{(k)}(u)|\leq 2 |Q^{(k)}(m)| \leq 2c_{k}^{-1} C m^{-\beta /2} N_{k}^{1-\frac{\beta}{2} - \alpha_k }.
\end{equation*}
To relate the quantities $u^{-\beta /2}$ and $m^{-\beta /2}$, we use the same method applied to the function $x^{-\beta /2}$. For $u\in [m, m+1/2]$:
\begin{equation*}
m^{-\beta/2}- u^{-\beta/2}\leq \frac{1}{2}\max_{x\in [m,m+1/2]}\left| \frac{d}{dx}x^{-\beta/2}\right|,
\end{equation*}
and therefore
\begin{equation*}
m^{-\beta/2}\leq u^{-\beta/2}+ \frac{\beta}{4}m^{-3\beta/2}.
\end{equation*}
The simple bounds $m^3 >  u$ and $\beta /4<1$ imply that
\begin{equation*}
\frac{\beta}{4} m^{-3\beta/2}\leq \frac{\beta}{4} u^{-\beta/2}\leq u^{-\beta/2}, 
\end{equation*}  
and thus
\begin{equation*}
|Q^{(k)}(u)| \leq 4c_{k}^{-1}Cu^{-\beta /2}N_{k}^{1-\frac{\beta}{2}-\alpha_k}, \quad u\in [m,m+1/2].
\end{equation*}
Having established this for $u\in [m, m+1/2]$, we can repeat the method on the interval $[m+1/2, m+1)$. By adjusting the constant $C$, the inequality \eqref{eq2.8} now holds with $m$ replaced by the continuous variable $u\geq 1$.

 In order to estimate $\hat{\mu}$, notice that for any $p\geq 1$,
\begin{equation}\label{eq2.9}
|\hat{\mu}(u)| \leq \prod_{k=1}^{p}|Q^{(k+1)}(\zeta_1 \cdots \zeta_k u)|.
\end{equation}

Given some $p$, we first consider values of $u$ of the form
\begin{equation}\label{eq2.9.1}
2^{p} N_1 \cdots N_{p} N_{p} \leq u <2^{p}N_1 \cdots N_{p}N_{p+1}.
\end{equation}
This corresponds to 
\begin{equation*}
N_{p} \leq \zeta_1 \cdots \zeta_{p} u < N_{p+1}.
\end{equation*}

Since the first $p-1$ terms of \eqref{eq2.9} may be bounded by $1$, it follows that 
\begin{equation*}
|\hat{\mu}(u)| \leq  Cc_{p+1}^{-1}(\zeta_1 \cdots \zeta_{p}u)^{-\beta/2} N_{p+1}^{1-\frac{\beta}{2}-\alpha_{p+1}}.
\end{equation*}
(By choosing $u$ in this range, we avoid the problematic cases of the argument of $|Q^{(p+1)}|$ being a multiple of $N_{p+1}$ or being very small, modulo $N_{p+1}$.) Because of the choice of the $\zeta_i$ in our construction, this gives
\begin{equation*}
|\hat{\mu}(u)| \leq c_{p+1}^{-1}C 2^{p\beta /2} N_{1}^{\beta /2}\cdots N_{p}^{\beta/2 } u^{-\beta /2}N_{p+1}^{1-\frac{\beta}{2} - \alpha_{p+1}}.
\end{equation*}

As stated in Section \ref{sec2.2}, we may assume that $c_k > \sqrt{2}C$ for each $k\in \mathbb{N}$, which means that 
\begin{equation}\label{eq2.10}
|\hat{\mu}(u)|\leq 2^{p\beta /2}u^{-\beta /2}N_{1}^{\beta /2}\cdots N_{p}^{\beta /2} N_{p+1}^{1-\frac{\beta}{2} - \alpha_{p+1}}.
\end{equation}
By requiring the terms $N_k$ to rapidly increase (passing to a subsequence if necessary), we can assume that 
\begin{equation*}
2^{p\beta /2}N^{\beta/2}_{1}  \cdots N^{\beta/2}_{p} N_{p+1}^{1-\frac{\beta }{2} - \alpha_{p+1}} = o(1),
\end{equation*}
This is possible, since by assumption $\beta > 2-2\alpha$. For given $\beta$ and $\alpha_i \to \alpha $ in the prescribed manner, we also have that $\beta >2-2\alpha_i$ for $\alpha_i$ close enough to $\alpha$, which we assume is the case, and hence that $1-\beta /2 -\alpha_{p+1} <0$. 

Now consider $u$ such that
\begin{equation*}
2^{p}N_1 \cdots N_{p}N_{p+1} \leq u < 2^{p+1} N_1 \cdots N_{p} N_{p+1}^2,
\end{equation*}
which corresponds to
\begin{equation*}
\frac{1}{2} \leq \zeta_1 \cdots \zeta_{p+1} u < N_{p+1}.
\end{equation*}
For this range of $u$ we apply the above argument, only with the approximation of $|Q^{(p+2)}|$ instead of $|Q^{(p+1)}|$. The case
$1/2 \leq  \zeta_1 \cdots \zeta_{p+1} u < 1$ should be mentioned, as the previous calculation required that $\zeta_1 \cdots \zeta_{p+1} u \geq 1$. However, since the upper bound on the derivative is valid for all $u>0$, we can similarly approximate $|Q^{(k)}|$:
\begin{equation*}
\max_{u\in [1/2,1)}|Q^{(k)}(u)| \leq 2 |Q^{(k)}(1)| \leq 2c_{k}^{-1}C N_{k}^{1-\frac{\beta}{2}-\alpha_k} .
\end{equation*}
Since $u^{-\beta/2}>1$ for $u \in [1/2,1)$, we see that 
\begin{equation*}
|Q^{(k)}(u)| \leq 2c_{k}^{-1}C u^{-\beta /2} N_{k}^{1-\frac{\beta}{2}-\alpha_k}
\end{equation*}
on this interval. We therefore obtain 
\begin{equation*}
|\hat{\mu}(u)|\leq 2^{(p+1)\beta /2}u^{-\beta /2}N_{1}^{\beta /2}\cdots N_{p+1}^{\beta /2} N_{p+2}^{1-\frac{\beta}{2} - \alpha_{p+2}}
\end{equation*}
for all $u$ in the current range. The required decay is ensured by the rapid increase of $N_k$.

By using increasing values of $p$, we see that 
\begin{equation*}
|\hat{\mu}(u)| = o(u^{-\beta /2})
\end{equation*}
for all $\beta $ satisfying the conditions of Theorem \ref{thm2.2}. 
\end{proof}

The result indicates that at least some sets satisfying the conditions of Theorem \ref{thm2.2} will be mapped to sets satisfying the conditions of Theorem \ref{thm2.1}, meaning that both the original set and its image will contain 3-arithmetic progressions. The next section further explores the invariance of the existence of arithmetic progressions under the mapping.

\subsection{Preservation of arithmetic progressions}\label{sec2.4}
The purpose of this section is to show that arithmetic progressions are preserved when transitioning from the integers to the continuum as described in previous sections. This is not dependent on the Hausdorff- or Fourier-dimensional aspects of the sets, only on the construction. We shall frequently refer to an arithmetic progression of length $n$ as an $n$AP.

Although the set $B(A)$ is not uniquely determined by $A$ alone but differs depending on constants used throughout, we suppress mention of the constants and regard $B(A)$ as fixed, and refer merely to $B$ throughout the proof.

\begin{proof}[Proof of Proposition 1.4]
 We first show that 3APs occurring in a subset of $[0,1]$ according to the construction \eqref{eq2.4} must arise from a 3AP in the integers. Recall that every element of our set $B$ can be written in the form \eqref{eq2.6}. Suppose that there are $a_1, a_2, a_3 \in B$ such that $a_2 - a_1 = a_3 - a_2$. Each of the $a_i$, $i=1,2,3$, is the image of a sequence of integers under the construction of $B$, and we can write
\begin{equation*}
a_ i =   \frac{a^{(i)}_1}{N_1 } + \frac{\zeta_1 a^{(i)}_2}{N_2 }+\frac{\zeta_1 \zeta_2  a^{(i)}_3}{N_3}+\cdots ,
\end{equation*}
where $a^{(i)}_j \in [1,N_i)\cap A$ for $i=1,2,3$; $j\in \mathbb{N}$.  

We therefore have that
\begin{equation*}
\frac{2a^{(2)}_1 - a^{(1)}_1 - a^{(3)}_1}{N_1 } + \frac{\zeta_1 (2a^{(2)}_2 - a^{(1)}_2 - a^{(3)}_2)}{N_2 }+\frac{\zeta_1 \zeta_2  (2a^{(2)}_3 - a^{(1)}_3 - a^{(3)}_3)}{N_3}+\cdots = 0.
\end{equation*}

Let $k\in \mathbb{N}$ be such that the 3AP occurs within an interval of length $\zeta_1 \cdots \zeta_{k-1}$ but not within an interval of length $\zeta_1 \cdots \zeta_{k}$. The choices of the parameters $\zeta_i$ guarantee that such a $k$ must exist, since the largest distance between any two points in a new subdivision of an interval is less than the distance to the next interval.  This would imply that the first $k-1$ terms of the expression above will be $0$ (trivially, since the points have equal terms in the expansion up to that point) and that the $k$-th term will only be $0$ should $2a_{k}^{(2)}-a_{k}^{(1)}-a_{k}^{(3)}=0$.
If the $k$-th term is not $0$, we have either that $a_{k}^{(3)} - a_{k}^{(2)} \geq a_{k}^{(2)} - a_{k}^{(1)} +1$ or that 
$a_{k}^{(3)} - a_{k}^{(2)} \leq a_{k}^{(2)} - a_{k}^{(1)} -1$. Suppose first that the former of the two alternatives holds, and consider the restricted case $a_{k}^{(3)} - a_{k}^{(2)} = a_{k}^{(2)} - a_{k}^{(1)} +1$.

The maximum distance between two points in the two separate intervals indexed by $a_{k}^{(1)}$ and $a_{k}^{(2)}$ is strictly less than 
\begin{equation*}
\zeta_1 \cdots \zeta_{k-1} \frac{a_{k}^{(2)} - a_{k}^{(1)}}{N_k } + \zeta_1 \cdots \zeta_k,
\end{equation*}
remembering that the intervals are half-open.
 The minimum distance between two points in the intervals indexed by $a_{k}^{(2)}$ and $a_{k}^{(3)}$ is strictly greater than 
\begin{align*}
& \zeta_1 \cdots \zeta_{k-1} \frac{ a_{k}^{(2)} - a_{k}^{(1)}+1}{N_k } - \zeta_1 \cdots \zeta_k \\ &= \zeta_1 \cdots \zeta_{k-1} \frac{a_{k}^{(2)} - a_{k}^{(1)}}{N_k } + \zeta_1 \cdots \zeta_{k-1 }\left( \frac{1}{N_k } - \zeta_k \right).  
\end{align*}
Requiring, as we did in the construction, that $\zeta_k = 1/2N_k $, we see that it is impossible for the difference between any two elements in the intervals indexed by $a_{k}^{(2)}$ and $a_{k}^{(1)}$ to be added to any element in the interval indexed by $a_{k}^{(2)}$ to obtain an element of the interval indexed by $a_{k}^{(3)}$. Therefore, it cannot be true that $a_{k}^{(3)} - a_{k}^{(2)} = a_{k}^{(2)} - a_{k}^{(1)} +1$. Since increasing $k$ only increases the minimum distance, it is not true that $a_{k}^{(3)} - a_{k}^{(2)} = a_{k}^{(2)} - a_{k}^{(1)} +k$ for any integer $k\geq 1$.

The same reasoning discounts the case $a_{k}^{(3)} - a_{k}^{(2)} \leq a_{k}^{(2)} - a_{k}^{(1)}-1 $, by considering the minimum possible  distance between points in the two separate intervals indexed by $a_{k}^{(2)}$ and $a_{k}^{(1)}$ and the maximum possible distance between points in the intervals indexed by $a_{k}^{(2)}$ and $a_{k}^{(3)}$. The minimum will always exceed the maximum, as in the previous case. Hence, $a_{k}^{(3)} - a_{k}^{(2)} = a_{k}^{(2)} - a_{k}^{(1)}$. Therefore, we can conclude that a 3AP in the set $B$ must be the result of a 3AP in the generating set $A$.

To pass to the general case of a $k$AP, we only need notice that such an AP can be seen as $k-2$ 3APs, each with the same common difference and the initial point of each shifted by the common difference each time. Here too, it is important to note the distances between intervals determined by the $\zeta_i$. Thus, the second 3AP must also correspond to a 3AP in $A$, with initial point the second member of the first 3AP, and so on, yielding a $k$AP in $A$. 

We now turn to the reverse, to show that a $k$AP in the generating set always leads to a $k$AP in $B$. Suppose that $P = \{ p_1 ,\ldots ,p_k \}\subseteq A $, $p_1 < p_2 < \ldots < p_k$, is a  $k$AP with common difference $d$ occurring in the $m$-th stage of the construction, that is, $p_k < N_m$. Then, for any fixed $a_{i} \in A\cap [0, N_i )$, $i\in \mathbb{N}\setminus \{ m \}$ there will be elements of $B$ expressible as 
\begin{equation*}
c_n =   \frac{a_1}{N_1 } + \zeta_1 \frac{ a_2}{N_2 }+\cdots + \zeta_1 \cdots \zeta_{m-1} \frac{p_n}{N_m} + \zeta_1 \cdots \zeta_m \frac{a_{m+1}}{N^{m+1}} + \cdots
\end{equation*} 
for $n=1,2,\ldots ,k$. Then $c_1 , \ldots ,c_k $ form a $k$AP in $[0,1]$ with common difference $\zeta_1 \cdots \zeta_{m-1}d/N^{m}$. 
\end{proof}

\section{Correspondence in higher dimensions}

\subsection{Hausdorff and Fourier dimensions}\label{sec3.1}

The method for obtaining subsets of $[0,1]^n$ from subsets of $\mathbb{N}^n$ is not much different from the one-dimensional case. First, we need to define the higher-dimensional notion of fractional density:
\begin{defn}\label{def3.1}
We say that a set $A\subseteq \mathbb{N}^n$ has upper fractional density $\alpha$, $0 \leq \alpha \leq n$, if
\begin{equation*}
\limsup_{N \to \infty} \frac{|A\cap [0,N )\times \cdots \times [0,N )|}{N^{\beta }}
\end{equation*}
is $\infty$ for any $\beta < \alpha$ and $0$ for any $\beta > \alpha$. This is indicated by $d_{f}^{*}(A) = \alpha$.
\end{defn}
The notions of \emph{lower fractional} density and \emph{fractional density} are defined analogously to the one-dimensional case.

With $A$ given, let $A_N = A\cap [0,N )\times \cdots \times [0,N )$, for some $N\in \mathbb{N}$. If a sequence of natural numbers $(N_i)_{i\geq 1}$ is given, we will set $A_i = A\cap [0,N_i )\times \cdots \times [0,N_i )$. The definition of a Salem-type set in $\mathbb{N}^n$ is as follows:
\begin{defn}\label{def3.2}
A set $A\subseteq \mathbb{N}^n$ is a Salem-type set if it has upper fractional density $\alpha \in (0,n]$ and for each $\vec{m} = (m_1 , \ldots ,m_n )\in [0,N)^n \setminus \{ \vec{0}\}$,
\begin{equation*}
\frac{1}{N^{n}} \sum_{\vec{n}\in [0,N )\times \cdots \times [0,N) } \chi_{A_N} (\vec{n}) e^{-2\pi i \frac{\vec{m}\cdot \vec{n} }{N}}  \leq C (|\vec{m}|N )^{-\beta/2}
\end{equation*}
for all large $N\in \mathbb{N}$, for some $C>0$ and $0< \beta \leq \alpha$. The supremum of the $\beta$ for which the above holds is the \emph{Fourier dimension} of $A$.
\end{defn}

The fact that the Fourier dimension will always be less than the fractional density follows from the correspondences in the next section. 

\subsection{Constructing a subset of $[0,1]^n$}\label{sec3.2}

Suppose that $A\subseteq \mathbb{Z}^n$ such that $d_{f}(A) = \alpha >0$. As before, we choose a strictly increasing sequence $( \alpha_i )_{i\geq 1}$ of real numbers such that $\alpha_i \to \alpha $, and a strictly increasing sequence $(\gamma_i )_{i\geq 1}$ tending to infinity. Let $(N_{i})_{i\geq 1}$ be an increasing sequence for which  
\begin{equation*}
c_i := \frac{|A\cap [0,N_{i} )\times \cdots \times [0,N_{i} )|}{N_{i}^{\alpha_i}} > \gamma_i.
\end{equation*} 
For each $i \in \mathbb{N}$, choose  $\zeta_{i} = (2N_i)^{-1}$ for all $i\in \mathbb{N}$. Define $A_i = A \cap [0,N_{i} ) \times \cdots \times [0, N_{i}) \subset \mathbb{N}^n$ for all $i\in \mathbb{N}$. Set
\begin{equation*}
B_1 = \bigcup_{(i_1 , \ldots ,i_n )\in A_1} \left[ \frac{i_1}{N_{1}}, \frac{i_1}{N_{1}} + \zeta_{1} \right) \times \cdots \times  \left[ \frac{i_n}{N_{1}}, \frac{i_n}{N_{1}} + \zeta_{1} \right).
\end{equation*}  
Define the set $J_k$  by 
\begin{equation*}
J_k = \left\{ (x_1, \ldots, x_n) \in \mathbb{R}^n: [x_1 ,x_1 +\zeta_1 \cdots \zeta_k) \times \cdots \times  [ x_n, x_n + \zeta_1 \cdots \zeta_k) \subseteq B_k \right\}. 
\end{equation*}
Once again, we define $B_{k+1}$ recursively:
\begin{equation*}
\begin{split}
B_{k+1} = \bigcup_{\vec{x}\in J_k } \bigcup_{\vec{i}\in A_{k+1} } 
\left[ x_1  +\frac{\zeta_1 \cdots \zeta_k i_1 }{N_{k+1}}, x_1  +\frac{\zeta_1 \cdots \zeta_k i_1 }{N_{k+1}} + \zeta_1 \cdots \zeta_{k+1}\right) \times \cdots
\times \\ \left[ x_n  +\frac{\zeta_1 \cdots \zeta_k i_n }{N_{k+1}}, x_n  +\frac{\zeta_1 \cdots \zeta_k i_n }{ N_{k+1}} + \zeta_1 \cdots \zeta_{k+1}  \right).
\end{split}
\end{equation*}

Finally, let
\begin{equation*}
B = \bigcap_{k=1}^{\infty} B_k.
\end{equation*}
As before, each element of the set $B$ can be expressed as a series
\begin{equation*}
\frac{\vec{a}^{(1)} (\varepsilon_1 )}{N_1} + \frac{\zeta_1 \vec{a}^{(2)} (\varepsilon_2 )}{N_2} + \frac{\zeta_1 \zeta_2 \vec{a}^{(3)} (\varepsilon_3 )}{N_3} + \cdots,
\end{equation*}
where $\vec{a}^{(k)}(i) = (a^{(k)}_{1}(i), \ldots ,a^{(k)}_{n}(i))$ denotes the $i$-th element of the set $A_k$ (with some ordering) and each $\varepsilon_i $ ranges over all values  $1,2, \ldots ,|A_i|.$ 

\begin{prop}\label{prop3.1}
Suppose that $A\subseteq \mathbb{N}^n$ has upper fractional density $\alpha$ and satisfies the conditions of Definition 3.2 for some $\beta$, where $\beta > 2n - 2\alpha$. Then the set $B$ constructed above has Hausdorff dimension $\alpha$ and Fourier dimension no less than $\beta$. 
\end{prop}

\begin{proof} Verifying that the set $B$ has Hausdorff dimension equal to the upper fractional density is no harder than in the one-dimensional case. For the existence of the measure and the decay of the Fourier-Stieltjes transform, we proceed analogously to the one-dimensional case, with only minor differences in the approximation of $\hat{\mu}$.  At each stage $k$ of the construction, a probability measure $\mu_k$ can be defined on the set $B_k$ which assigns equal measure to each cube comprising $B_k$, with a distribution function
\begin{equation*}
F_k (\vec{x}) = \mu_k ([\vec{y}\in [0,1]^n : y_i \leq x_i, i=1,\ldots ,n]).
\end{equation*} 
The $\mu_k$ converge weakly to a measure $\mu$. Setting 
\begin{equation*}
Q^{(k)}(\vec{u} ) = \frac{1}{|A_k |}\sum_{j=1}^{|A_k |} e^{-\frac{2\pi i}{N_{k}^{}} \vec{u} \cdot \vec{a}^{(k)}( j )},
\end{equation*}
we can write 
\begin{equation*}
\hat{\mu }(\vec{u}) = \lim_{n\to \infty} Q^{(1)}(\vec{u})\prod_{k=1}^{m} Q^{(k+1)} (\zeta_1 \cdots \zeta_k \vec{u})
\end{equation*}
and 
\begin{equation*}
| \hat{\mu }(\vec{u})| \leq \prod_{k=1}^{p}|Q^{(k+1)}(\zeta_1 \cdots \zeta_k \vec{u})|,
\end{equation*}
as before. We have the upper bound
\begin{equation*}
|Q^{(k)}(\vec{m}) | \leq \frac{C}{c_k} |\vec{m}|^{-\beta /2} N^{n-\frac{\beta}{2} - \alpha_k}
\end{equation*}
for $\vec{m} \in [1,N_k) \times \cdots \times [1,N_k) \subset \mathbb{N}^n$, and we can assume that this holds for non-integer arguments for much the same reasons as before. Specifically, for $\vec{x}=(x_1,\ldots ,x_n)\in \mathbb{R}^n$ and $\vec{u}=(u_1,\ldots ,u_n)$,
\begin{equation*}
\left| \frac{\partial}{\partial u_i}|Q^{(k)}(\vec{u})|\right| \leq |Q^{(k)}(\vec{u})|
\end{equation*} 
for each $i=1,\ldots ,n$. For a given $y\in \mathbb{R}$, let $\vec{u}_i (y) = (u_1,\ldots,u_{i-1}, y,u_{i+1} , \ldots ,u_n )$. The previous inequality implies that 
\begin{equation*}
\max_{y \in [u_i,u_i +a]}\left| \partial_{x_i}|Q^{(k)}(\vec{x})|\at[\big]{{\vec{x} = \vec{u}_i (y)} }\right|\leq \max_{y\in [u_i,u_i +a]} |Q^{(k)}(\vec{u}_i (y))| \leq \max_{\substack{0\leq  x_i \leq a \\ 1 \leq i \leq n}}|Q^{(k)}(\vec{u}+\vec{x})|,
\end{equation*}
and consequently
\begin{equation*}
\begin{aligned}
\max_{\substack{0\leq  x_i \leq a \\ 1 \leq i \leq n}}|Q^{(k)}(\vec{u}+\vec{x})| &\leq a \left( \sum_{i=1}^{n} \left( \max_{y \in [u_i,u_i +a]}| \partial_{x_i}|Q^{(k)}(\vec{x})|\at[\big]{{\vec{x} = \vec{u}_i (y)} } |_{} \right)^2 \right)^{1/2}+ |Q^{(k)}(\vec{u})|\\
& \leq  a \sqrt{n} \max_{\substack{0\leq  x_i \leq a \\ 1 \leq i \leq n}}|Q^{(k)}(\vec{u}+\vec{x})| + |Q^{(k)}(\vec{u})|,
\end{aligned}
\end{equation*}
yielding
\begin{equation*}
\max_{\substack{0\leq  x_i \leq a \\ 1 \leq i \leq n}}|Q^{(k)}(\vec{u}+\vec{x})| \leq \frac{1}{1-a\sqrt{n}} |Q^{(k)}(\vec{u})|
\end{equation*}
as long as $1-a\sqrt{n}>0$. As with the one-dimensional case, we can now repeatedly extend the upper bound over cubes of side-length no greater than, for instance, $(2\sqrt{n})^{-1}$, by using the inequality 
\begin{equation*}
| \vec{m}|^{-\beta /2} \leq K|\vec{u}|^{-\beta /2}
\end{equation*}
 for some constant $K$ independent of $\vec{m}$ and $\vec{u}$, with $\vec{u}$ contained in the cube
\begin{equation*}
\{ \vec{u} = (u_1, \ldots ,u_n) \in \mathbb{R}^n: 0 \leq u_i - m_i \leq (2\sqrt{n})^{-1},  i=1,2,\ldots ,n \}.
\end{equation*}
With the upper bound on $|Q^{(k)}(\vec{u})|$ established for this cube, we can repeat the argument to obtain the same bound for neighboring cubes, and therefore for all $\vec{u} \in \mathbb{R}^n$ with $u_i \geq 1$, $i=1,\ldots, n$. When $ u_i \in [1/2, 1)$, the argument in the one-dimensional case applies similarly.

The reasoning of Section \ref{sec2.3} now transfers to the multidimensional case almost seamlessly. The restrictions on the values of $u$ in the previous must only be replaced by restrictions on the coordinates of $\vec{u}$; thus, we consider the cases
\begin{equation*}
2^{p} N_1 \cdots N_{p} N_{p} \leq u_i <2^{p}N_1 \cdots N_{p}N_{p+1}, \quad i=1,2,\ldots ,n
\end{equation*}
and
\begin{equation*}
2^{p} N_1 \cdots N_{p} N_{p+1} \leq u_i <2^{p+1}N_1 \cdots N_{p}N_{p+1}^{2}, \quad i=1,2,\ldots ,n.
\end{equation*}
Also taking into account that now we have $\beta > 2n-2\alpha$, we find that  
\begin{equation*}
|\hat{\mu}(\vec{u})| = o(|\vec{u}|^{-\beta /2}).
\end{equation*}
We can conclude that the Fourier dimension of $B$ is $\beta$ or greater. 
\end{proof}

\section{Finite configurations in sparse subsets}

\subsection{Subsets of $\mathbb{R}^n$}\label{sec4.1}

In this section we briefly describe the result of Chan, \L{}aba and Pramanik~\cite{chan2016finite}, to which we shall find a corresponding theorem in the integers. Whereas before we used vector notation to indicate elements of higher-dimensional spaces in order to distinguish from the one-dimensional case, we abandon that now for readability.

\begin{defn}\label{def4.1}~\cite{chan2016finite}
Fix integers $n\geq 2$, $k\geq 3$ and $m\geq n$. Suppose $\mathbb{B} = \{B_1 ,\ldots ,B_k\}$ is a collection of $n\times (m-n)$ real matrices.
\begin{itemize}
\item[(a)] We say that $E\subset \mathbb{R}^n$ contains a $\mathbb{B}$-configuration if there exist $x\in \mathbb{R}^n$ and $y\in \mathbb{R}^{m-n}\setminus \{0\}$ such that $\{ x+B_j y \}_{j=1}^{k} \subseteq E$. 
\item[(b)] Given any finite collection of subspaces $V_1 ,\ldots ,V_q \subseteq \mathbb{R}^{m-n}$ with $\textrm{dim} (V_i ) < m-n$, we say that $E$ contains a non-trivial $\mathbb{B}$-configuration with respect to $(V_1 ,\ldots ,V_q )$ if there exist $x\in \mathbb{
R}^n $ and $y\in \mathbb{R}^{m-n}\setminus \bigcup_{i=1}^{q}V_i $ such that $\{ x+B_j y\}_{j=1}^{k} \subseteq E$.
\end{itemize}
\end{defn}

In order to formulate the theorem, we shall need the notion of \emph{non-degeneracy} for our configurations. Let $A_1 ,\ldots ,A_k$ be $n\times m$ matrices. For any set of (distinct) indices $J=\{j_1 ,\ldots ,j_s \} \subseteq \{ 1,\ldots ,k\}$, define the $ns\times m$ matrix $\mathbb{A}_J$ by
\begin{equation*}
\mathbb{A}_{J}^{t} = (A_{j_1}^t \cdots A_{j_s}^{t}).
\end{equation*}

Let $r$ be the unique positive integer such that 
\begin{equation*}
n(r-1) < nk-m \leq nr.
\end{equation*}

The meaning of this integer is as follows. Suppose that we construct an $m\times m$ matrix from certain of the matrices $A_{1}^{t} ,\ldots ,A_{k}^{t}$, each of which consists of $n$ columns. If $m$ is not a multiple of $n$, we fill up the remaining columns (fewer than $n$) with columns from one of the other matrices. The total number of columns of all matrices minus the number used for the $m\times m$ matrix is therefore $nk - m$. The number $r-1$ then represents the maximum number of groups of $n$ columns that can be entirely contained in $nk-m$ columns, unless $nk-m$ is a multiple of $n$, in which case $r$ groups are contained. In other words, $r$ is the total number of groups of $n$ columns that are wholly or partly contained in the $nk-m$ columns not contained in the $m\times m$ matrix. The number of columns fewer than $n$, outside these $r$ groups, that remain ``outside" the $m\times m$ matrix constructed, is then
\begin{equation*}
n' = nk - m - (r-1)n.
\end{equation*}   
It is possible that $n' = n$, when $nk-m$ is a multiple of $n$. 

\begin{defn}\label{def4.2}\cite{chan2016finite}
We say that $\{ A_1 ,\ldots ,A_k \}$ is \emph{non-degenerate} if for any $J\subseteq \{1,\ldots ,k\} $ with $|J| = k-r$ and any $j\in \{1,\ldots ,k\}\setminus J$, the $m\times m$ matrix
\begin{equation*}
(\mathbb{A}_{J}^{t} \tilde{A_j}^{t} )
\end{equation*}
is non-singular for any choice $\tilde{A_j}^{t}$ a submatrix of $A_j$ consisting of $n-n'$ rows.
\end{defn}
That is, we create an $m\times m$ matrix by using as many of the matrices $A_{i}^{t}$ as we can without obtaining more than $m$ columns in total, then fill out the rest of the matrix by using columns from one of the matrices not yet used. If all such matrices are non-singular, the assembly $\{ A_1 ,\ldots ,A_k \}$ is non-degenerate.

We then have the following theorem from Chan et al.\ \cite{chan2016finite}. Note that $I_{n\times n}$ denotes the $n\times n$ identity matrix.
\begin{thm}\label{thm4.1}
Suppose
\begin{equation*}
n \left\lceil \frac{k+1}{2} \right\rceil \leq m \leq nk
\end{equation*}
and 
\begin{equation*}
\frac{2(nk-m)}{k} < \beta <n.
\end{equation*}
Let $\mathbb{B} = \{B_1 ,\ldots ,B_k \}$ be a collection of $n\times (m-n)$ matrices such that $A_j = (I_{n\times n} B_j )$ is non-degenerate in the sense 
of Definition \ref{def4.2}. Then for any constant $C>0$, there exists a positive number $\varepsilon = \varepsilon (C,n,k,m,\mathbb{B} ) \ll 1$ with the following property. Suppose the set $E\subseteq \mathbb{R}^n $ with $|E| = 0$ supports a positive, finite, Radon measure $\mu $ satisfying the two conditions:
\begin{itemize}
\item[(a)]{ball condition:
\begin{equation*}
\sup_{\substack{x\in E \\ 0<r<1}} \frac{\mu(B(x,r))}{r^{\alpha}} \leq C \text{ if } n-\varepsilon < \alpha <n,
\end{equation*}}
\item[(b)]{Fourier decay:
\begin{equation*}
| \hat{\mu } (\xi) | \leq \frac{C}{(1+|\xi |)^{\beta /2}}.
\end{equation*} }
for every $\xi \in \mathbb{R}^n$.
\end{itemize}
Then
\begin{itemize}
\item[(i)]{E contains a $\mathbb{B}$-configuration.}
\item[(ii)]{For any finite collection of subspaces $V_1 , \ldots ,V_q \subseteq \mathbb{R}^{m-n}$ with $\dim (V_i ) < m-n$, $E$ contains a non-trivial $\mathbb{B} $-configuration with respect to $(V_1 ,\ldots ,V_q )$ in the sense of Definition \ref{def4.1}. }
\end{itemize}
 
\end{thm}

We remark that condition (b) is fulfilled when $|\hat{\mu }(\xi) | = o(|\xi |^{-\beta /2})$.

\subsection{Subsets of $\mathbb{Z}^n$}\label{sec4.2}

Matrices of the form specified in Theorem \ref{thm4.2} will be said to define ``par\-al\-lel\-o\-gram-like" configurations. This will be further illustrated in Section 4.3, which will also discuss the issue of degeneracy in the configurations, such as when some of the points of the parallelogram are equal, leading to a trivial configuration.

\begin{proof}[Proof of Theorem 1.5.]
 The proof is already largely accomplished, due to the constructions and correspondences of previous sections. Suppose that $A$ satisfies the conditions of the theorem. Let $B$ be the subset of $[0,1]^n$ obtained via the construction in Section \ref{sec3.2}.
By the dimensional correspondence of the construction, we see immediately that $B$ satisfies the conditions of Theorem \ref{thm4.1}, and must contain a $\mathbb{B}$-configuration. The rest of the proof aims to show that this configuration corresponds to one in $A \subseteq \mathbb{N}^n$, much as we earlier showed that arithmetic progressions in $B\subseteq [0,1]$ must correspond to arithmetic progressions in the generating set $A\subseteq \mathbb{N}$.

Without loss of generality, we suppose that $N_1$ is large enough to separate the points of the configuration (since no matter which $N_1$ we start with, $B$ must contain the required configuration). Partition $\mathbb{R}^n$ and $\mathbb{R}^{m-n}$ into cubes so that 
\begin{equation*}
\mathbb{R}^l = \bigcup_{z \in \mathbb{Z}^l} \bigcup_{i_1 ,\ldots ,i_l=0}^{N_1 -1}z + \left[ \frac{i_1}{N_1}, \frac{i_1 +1}{N_1}\right) \times \cdots \times
\left[ \frac{i_l}{N_1}, \frac{i_l+1}{N_1}\right),
\end{equation*}
for $l=n ,m-n$.

Let $x=(x_1 ,\ldots ,x_n) \in \mathbb{R}^n$ and $y = (y_1, \ldots ,y_{m-n}) \in \mathbb{R}^{m-n}$ be the two points whose existence is guaranteed by Theorem \ref{thm4.1}. By assumption, $B_1$ is the zero matrix, and hence $x\in B$. The $N_1$-approximation $x^{(1)} = (x^{(1)}_1 , \ldots ,x^{(1)}_n)\in \mathbb{Z}^n$ of the point $x$ in $\mathbb{R}^n$ is then the corner of the cube with side-lengths $1/N_1$ containing the point with the smallest coordinates in each dimension (in two dimensions, this would be the southwest corner of the containing square).  More precisely,
\begin{equation*}
\text{for each } i=1,2,\ldots ,n:  0\leq x_i - x^{(1)}_{i} < \frac{1}{2N_1} \text{ and } \exists a_i \in \mathbb{N} \text{ such that } x^{(1)}_{i} = \frac{a_i}{N_1}.
\end{equation*}
Note that the nature of the construction of $B$, specifically the choice of $\zeta_1$, ensures that this approximation is possible and unique.
For any $z\in \mathbb{R}^n$, $z^{(1)}$ shall denote the $N_1$-approximation of $z$ in the above sense. Elements of $\mathbb{R}^{m-n}$ (considered as a distinct space for the purposes of the argument even when $n=m-n$), will be approximated slightly differently. Approximate $y\in \mathbb{R}^{m-n}$ by the closest corner available, that is, the point $y^{(2)}  = (y^{(2)}_1, \ldots , y^{(2)}_{m-n})$ such that 
\begin{equation*}
\text {for each } i=1,2,\ldots ,m-n:  |y_i - y^{(2)}_{i}| \leq \frac{1}{2N_1} \text{ and } \exists b_i \in \mathbb{N} \text{ such that } y^{(2)}_{i} = \frac{b_i}{N_1}.
\end{equation*} 
If, for any coordinate,  there are two possible choices for an approximation, we choose the one for which $y - y_{i}^{(2)} >0$, thus ensuring that $y -y_{i}^{(2)} >-(2N_1)^{-1}$ for all $i$.   

We can write 

\begin{equation*} x =
\begin{bmatrix}
x^{(1)}_1 + \epsilon(x^{(1)}_{1})\\
x^{(1)}_2 + \epsilon (x^{(1)}_{2})\\
\vdots \\
x^{(1)}_n + \epsilon (x^{(1)}_{n})
\end{bmatrix}, \quad
y = 
\begin{bmatrix}
y^{(2)}_1 + \epsilon(y^{(2)}_{1})\\
y^{(2)}_2 + \epsilon(y^{(2)}_{2})\\
\vdots \\
y^{(2)}_{m-n} + \epsilon(y^{(2)}_{m-n})
\end{bmatrix}, 
\end{equation*}
where  $ 0\leq \epsilon (x_{i}^{(1)}) <(2N_1)^{-1}$, $1\leq i \leq n$, and $-(2N_1)^{-1}< \epsilon ( y_{j}^{(2)}) \leq (2N_1)^{-1}$, $1\leq j \leq m-n$. We shall use the notation
\begin{align*}
\epsilon (x^{(1)}) &=   (\epsilon (x^{(1)}_{1}) ,\ldots ,\epsilon(x^{(1)}_{n}))^t \\
\epsilon (y^{(2)}) &=   (\epsilon (y^{(2)}_{1}) ,\ldots ,\epsilon(y^{(2)}_{m-n}))^t.  
\end{align*}

We assert that $x^{(1)}+B_j y^{(2)}$ is the $N_1$-approximation of the point $x+B_j y$, $j=1, \ldots ,k$.
Observe the $j$-th error term
\begin{equation*}
\epsilon(x^{(1)})+B_j \epsilon (y^{(2)}) =  \begin{bmatrix}
\epsilon({x^{(1)}_1})\\
\epsilon({x^{(1)}_2})\\
\vdots \\
\epsilon({x^{(1)}_n})
\end{bmatrix} + 
B_j \begin{bmatrix}
\epsilon({y^{(2)}_1})\\
\epsilon({y^{(2)}_2})\\
\vdots \\
\epsilon({y^{(2)}_{m-n}})
\end{bmatrix}.
\end{equation*}
Each component of this term represents the extent to which each coordinate of $x^{(1)}+B_j y^{(2)}$ differs from a multiple of $1/N_1$.
Our assertion is then that 
\begin{equation*}
0 \leq (\epsilon(x^{(1)})+B_j \epsilon (y^{(2)}))_i <\frac{1}{N_1}
\end{equation*}
for the $i$-th component of $\epsilon(x^{(1)})+B_j \epsilon (y^{(2)})$, $i=1,2,\ldots ,n$. For the moment, consider only $j=1,2,\ldots ,m/n$. 
For the assertion to be violated, we must have that 
\begin{equation*}
( \epsilon (x^{(1)}) + B_j \epsilon (y^{(2)}))_i \geq \frac{1}{N_1}
\end{equation*}
or
\begin{equation*}
( \epsilon (x^{(1)}) + B_j \epsilon (y^{(2)}))_i <0
\end{equation*}
for some value of $i$. The first case is eliminated by the bounds on the error terms and the restrictions on entries of $B_j$. Since we have $x+B_j y \in B$ (and must therefore be contained in one of the cubes specified by the construction at the first stage), the second case can only occur when
\begin{equation*}
( \epsilon (x^{(1)}) + B_j \epsilon (y^{(2)}))_i \leq -\frac{1}{2N_1}
\end{equation*}
in some coordinate $i$, which represents the case where the $N_1$-approximation $(x+B_j y )^{(1)}$ is less in some coordinate than $x^{(1)}+B_j y^{(1)}$. Once again, this is prohibited by the bounds on $\epsilon (x^{(1)})$ and $\epsilon (y^{(2)})$.

By construction, the point $x^{(1)}$ corresponds to a point in the original set $A$, specifically $x_A = (N_1 x^{(1)}_1 ,\ldots , N_1 x^{(1)}_n )$.
In turn, $y^{(2)}$ defines an element $y_A = (N_1 y^{(2)}_1 ,\ldots ,N_1 y^{(2)}_{m-n})$ of $\mathbb{Z}^{m-n}$. For each point $x^{(1)}+B_j y^{(2)}$, $j=1,2,\ldots ,m/n$, there therefore exists some $a_j\in A$ such that 
\begin{equation*}
x^{(1)}+B_j y^{(2)} = \frac{a_j}{N_1},
\end{equation*}
and we have that 
\begin{equation*}
x_A + B_j y_A = a_j , \quad j=1,2, \ldots ,m/n.
\end{equation*}
To deal with the matrices $B_j$, $m/n < j \leq k$,  observe again that the point $x+B_l y$, $l=1,\ldots ,m/n$, is an element of the set $B$. This implies that, for each coordinate $i$,   
\begin{equation*}
0 \leq ( \epsilon (x^{(1)}) + B_l \epsilon (y^{(2)}))_i < \frac{1}{2N_1}, \quad l= 1,2, \ldots ,m/n.
\end{equation*} 
Since for $j = m/n+1, \ldots ,k$ we have $B_j = B_{j_1} + B_{j_2}$ for some $j_1 ,j_2 \in \{2,3,\ldots ,m/n\}$, the above implies that 
\begin{equation*}
-\frac{1}{2N_1} < (\epsilon (x^{(1)}) + B_{j_1}\epsilon (y^{(2)}) + B_{j_2}\epsilon (y^{(2)}))_i < \frac{1}{N_1}
\end{equation*}
for any $i = 1,2,\ldots ,n$, and hence that
\begin{equation*}
-\frac{1}{2N_1}<( \epsilon (x^{(1)}) + B_j \epsilon (y^{(2)}))_i < \frac{1}{N_1}
\end{equation*} 
for any $j=m/n+1,\ldots ,k$. By the same argument as above, taking into account that $x+B_j y\in B$ also for $j= m/n +1, \ldots ,k$, we can conclude that
\begin{equation*}
0\leq ( \epsilon (x^{(1)}) + B_j \epsilon (y^{(2)}))_i < \frac{1}{N_1}
\end{equation*} 
and hence 
\begin{equation*}
x^{(1)}+B_j y^{(2)} = \frac{a_j}{N_1}
\end{equation*}
for some $a_j \in A$, $m/n < j \leq k$. This establishes the result. 
\end{proof}

\subsection{Existence of parallelograms}\label{sec4.3}

Fix $n\geq 2$, $k=4$ and $m=3n$ and suppose the set $A\subseteq \mathbb{Z}^n$ satisfies the conditions of Theorem \ref{thm4.2}. Let  
\begin{align*}
B_1 &= 0_{n\times 2n}\\
B_2 &= (I_{n\times n} 0_{n\times n})\\
B_3 &= (0_{n\times n} I_{n\times n})\\
B_4 &= B_2 + B_3.
\end{align*}
Chan et al.\ \cite{chan2016finite} show that the  collection of matrices $\{A_1 , A_2 ,A_3 ,A_4 \} $, formed from $\{ B_1 ,B_2 , B_3 ,B_4 \}$ in the manner described in Theorem \ref{thm4.2}, are non-degenerate and satisfy all other conditions of the theorem. This means that there are some $x\in \mathbb{Z}^n$, $y\in \mathbb{Z}^{2n}$ such that the points
\begin{equation*}
x, x+ \begin{bmatrix}
y_1 \\ \vdots \\ y_n
\end{bmatrix}, 
x+\begin{bmatrix}
y_{n+1} \\ \vdots \\ y_{2n}
\end{bmatrix},
x+\begin{bmatrix}
y_1 + y_{n+1} \\ \vdots \\ y_n + y_{2n }
\end{bmatrix}
\end{equation*}
are all contained in $A$. In other words, there are points $x, u, v \in \mathbb{Z}^n$ such that
\begin{equation*}
\{ x, x+u , x+v , x+u+v \} \subset A.
\end{equation*} 
It is important to mention the possibility of degenerate configurations, where not all points are distinct. Setting $v=0$ in the above configuration would mean that sets could trivially contain parallelograms as long as they consist of more than two elements. This is addressed by part $(ii)$ of the conclusion of Theorem \ref{thm4.1}. By defining the subspaces $V_1, V_2 , V_3 ,V_4$ suitably, we can conclude that the points obtained are all distinct; see Chan et al.\  \cite{chan2016finite}, Section $7$. Since the points are also separated by our $N_1$-approximation in Theorem \ref{thm4.2}, we can conclude that the elements of $A$ obtained are also distinct.

\subsection{Possible generalizations of configurations}\label{sec4.4}

The following generalizations of Theorem \ref{thm4.2} suggest themselves.

\begin{enumerate}
\item{There is of course no need to require that $y\in \mathbb{Z}^{m-n}$ in Theorem \ref{thm4.2}, nor even that the matrices must all have integer entries, as long as the configuration obtained exists in $\mathbb{Z}^n$. For instance, if we consider assemblies $\mathbb{B}$ in which the matrices have some rational but non-integer coefficients that satisfy the constraints of Theorem \ref{thm4.1}, such configurations could still exist in the set $B$ if conditions are satisfied. This would imply that we could consider $y\in \mathbb{Q}^{m-n}$ instead, for which the corresponding configuration exists in $A$. This could be useful in establishing analogues of Corollary 1.7 in Chan et al.\ \cite{chan2016finite}, concerning the existence of certain triangle configurations. }
\item{It is not necessary to assume that $\zeta_ i = (2N_i )^{-1}$ during the construction in Section 2. The $\zeta_i$ could be taken to be, for instance, $(3N_i )^{-1}$ (as long as some small extra growth condition is applied to $(N_i )$), which would open the door to proving Theorem \ref{thm4.2} for configurations with some matrices being the sum of three other matrices. This could allow one to prove the existence of parallelogram-like configurations in $A$ such as $\{x, x+u, x+v, x+z, x+u+v+z\}$, for instance.}
\item{Only a single value of $m$, the smallest possible, was considered in Theorem \ref{thm4.2}. Larger values should be considered for more general configurations, although the non-degeneracy conditions on the matrices $A_1 ,\ldots ,A_k$ complicate the proof of Theorem \ref{thm4.2}}
\item{Just as arithmetic progressions of length greater than $3$ seem to require higher-order uniformity conditions 
(see Conlon et al. \cite{conlon2015relative}, for instance), we can conjecture that more general configurations in $\mathbb{Z}^n$ and $[0,1]^n$ will require such conditions. It seems feasible, given the results of this paper, that these conditions will be preserved under the current construction.  }
\end{enumerate}

\section{Acknowledgements}

The author would like to thank an anonymous referee whose comments and suggestions greatly improved the presentation and correctness of this paper.

\bibliographystyle{plain}

\end{document}